\documentclass[a4paper,11pt,oneside]{article}
\usepackage{amsmath,amssymb,amsthm}
\usepackage{mathtools}

\DeclarePairedDelimiter\floor{\lfloor}{\rfloor}
\allowdisplaybreaks

\newcommand{\Z}{\mathbb{Z}}
\newcommand{\pres}[2]{\langle {#1}\ |\ {#2} \rangle}
\newcommand{\gpres}[1]{\langle {#1} \rangle}
\newcommand{\npres}[1]{\langle\langle {#1} \rangle\rangle}
\newtheorem{theorem}{Theorem}
\newtheorem{lemma}[theorem]{Lemma}
\newtheorem{example}[theorem]{Example}
\newtheorem{corollary}[theorem]{Corollary}
\numberwithin{theorem}{section}

\begin{document}
\title{Fractional Fibonacci groups with an odd number of generators}
\author{Ihechukwu Chinyere and Gerald Williams\thanks{This work was supported by the Leverhulme Trust Research Project Grant RPG-2017-334.}}

\maketitle

\begin{abstract}
The Fibonacci groups $F(n)$ are known to exhibit significantly different behaviour depending on the parity of $n$. We extend known results for $F(n)$ for odd $n$ to the family of Fractional Fibonacci groups $F^{k/l}(n)$. We show that for odd $n$ the group $F^{k/l}(n)$ is not the fundamental group of an orientable hyperbolic 3-orbifold of finite volume. We obtain results concerning the existence of torsion in the groups $F^{k/l}(n)$ (where $n$ is odd) paying particular attention to the groups $F^k(n)$ and $F^{k/l}(3)$, and observe consequences concerning the asphericity of relative presentations of their shift extensions. We show that if $F^{k}(n)$ (where $n$  is odd) and $F^{k/l}(3)$ are non-cyclic 3-manifold groups then they are isomorphic to the direct product of the quaternion group $Q_8$ and a finite cyclic group.
\end{abstract}

\medskip


\noindent \textbf{MSCs:} 20F05, 57M05, 57M60.\\

\section{Introduction}\label{sec:introduction}

The \em Fibonacci groups \em
\[ F(n) =\pres{x_0, \ldots ,x_{n-1}}{ x_ix_{i+1}=x_{i+2}\ (0\leq i<n)}\]
(subscripts mod~$n$) were introduced by Conway in \cite{Conway65}, and they have since been studied from both algebraic and topological perspectives. The \em Fractional Fibonacci groups \em
\begin{alignat}{1}
F^{k/l}(n) =\pres{x_0, \ldots ,x_{n-1}}{ x_i^lx_{i+1}^k=x_{i+2}^l\ (0\leq i<n)}\label{eq:Fkl(n)}
\end{alignat}
where $k,l\neq 0$, $n\geq 1$, subscripts mod~$n$, introduced in \cite{KimVesninSiberian}, generalise the Fibonacci groups $F(n)=F^{1/1}(n)$ and also the groups $F^{k}(n)=F^{k/1}(n)$ considered in \cite{Maclachlan,MaclachlanReid}. For even $n \geq 6$ and coprime integers $k,l\geq 1$ the groups $F^{k/l}(n)$ have been shown to be fundamental groups of 3-manifolds (see \cite{HKM,HLM,HLM2,CS} for the case $k=l=1$, see \cite{MaclachlanReid} for the case $l=1$, and \cite{KimVesninSiberian} for the case of coprime integers $k,l\geq 1$).

It is known that the Fibonacci groups $F(n)$ exhibit substantially different behaviour depending on the parity of $n$. For instance, if $n$ is even then $F(n)$ is the fundamental group of a 3-manifold, namely an $n/2$-fold cyclic cover of $S^3$ branched over the figure eight knot (which is spherical if $n=2,4$, an affine Riemannian manifold if $n=6$, and hyperbolic if $n\geq 8$) \cite{HKM,HLM,HLM2,CS}, whereas if $n\geq 3$ is odd then $F(n)$ is not the fundamental group of any hyperbolic 3-orbifold of finite volume \cite[Theorem~3.1]{Maclachlan}, and if $n\geq 9$ is odd then $F(n)$ is not the fundamental group of any 3-manifold \cite[Theorem~3]{HowieWilliamsMFD}. Moreover, if $n\geq 6$ is even then $F(n)$ is infinite and torsion-free by statements P(3),P(4) of \cite{HKM} whereas if $n\geq 9$ is odd then $F(n)$ is infinite \cite{HoltFibonacci,Newman,Lyndonunpublished,Chalk} and contains an element of order 2 by \cite[Proposition~3.1]{BardakovVesnin} ($F(2), F(3), F(4), F(5), F(7)$ are finite groups).

In this article we consider the Fractional Fibonacci groups $F^{k/l}(n)$ when $n$ is odd. In Section~\ref{sec:prelims} we obtain some basic observations about the groups $F^{k/l}(n)$. In Section~\ref{sec:abelianisation} we obtain a recurrence relation formula for the order $|F^{k/l}(n)^\mathrm{ab}|$ (Theorem~\ref{thm:Fk/labVn}) and consequences of it that will be used in later sections. In Section~\ref{sec:hyporbifolds} we prove Theorem~\ref{thm:hyper}, which states that for odd $n$ the group $F^{k/l}(n)$ is not the fundamental group of an orientable hyperbolic 3-orbifold of finite volume and in Corollary~\ref{cor:hyperkodd}, we prove that if, in addition, $k$ is odd then $F^{k/l}(n)$ is not the fundamental group of a hyperbolic 3-orbifold of finite volume. In Section~\ref{sec:torsion} we consider torsion elements in $F^{k/l}(n)$ and introduce a word $w(n,k)$ that is the basis for much of this section. By a result of Bardakov and Vesnin \cite{BardakovVesnin}, for odd $n\geq 9$, the word $w(n,1)$ is an element of order~2 in $F(n)$ (Theorem~\ref{thm:Worder2inF(n)}) and this has consequences for the asphericity of the relative presentation of the shift extension of $F(n)$ (Corollary~\ref{cor:F(n)extnotasp}), and the result was the basis for the proof in \cite{HowieWilliamsMFD} that $F(n)$ is not the fundamental group of a 3-manifold (Theorem~\ref{thm:F(n)not3mfd}). We develop extensions of these results to the general case $F^{k/l}(n)$ and apply them to the groups $F^k(n)$ and $F^{k/l}(3)$. In Theorem~\ref{thm:order2inFk/l} we show that $w(n,k)^2=1$ in $F^{k/l}(n)$ and that $w(n,k)$ is a commutator. Corollary~\ref{cor:Fkabelian} shows that $w(n,k)=1$ if and only if $F^{k}(n)$ is abelian, and Corollary~\ref{cor:Fk/l(3)abelian} does the same for the group $F^{k/l}(3)$. Corollary~\ref{cor:Fkl(n)extnotasp} then shows that if $w(n,k)\neq 1$ then the relative presentation of the shift extension of $F^{k/l}(n)$ is not aspherical and Corollaries~\ref{cor:Fk(n)extnotasp},\ref{cor:Fk/l(3)extnotasp} show that this relative presentation is not aspherical in the cases $l=1$ and $n=3$, respectively. In Section~\ref{sec:3mfd} we consider when $F^{k/l}(n)$ is a 3-manifold group and show that if $w(n,k)\neq 1$ then $F^{k/l}(n)$ is not a 2-generator, infinite, 3-manifold group (Lemma~\ref{lem:inf3mfd}). In Theorem~\ref{thm:Fkl(3)Fk(n)3mfd} we use this to prove that if $F^{k}(n)$ or $F^{k/l}(3)$ is a 3-manifold group then it is either a finite cyclic group or isomorphic to the direct product of the quaternion group $Q_8$ and a finite cyclic group.

\section{Preliminaries}\label{sec:prelims}

Our first lemma is immediate from the definition of $F^{k/l}(n)$.

\begin{lemma}\label{lem:easyFklgroups}
\begin{itemize}
  \item[(a)] $F^{k/l}(2)\cong \Z_k*\Z_k$;
  \item[(b)] $F^{k/0}(n)$ is isomorphic to the free product of $n$ copies of $\Z_k$.
\end{itemize}
\end{lemma}
For $n\geq 1$ and $l \in \Z$ let
\[ G(n,l)=\pres{x_0,\ldots,x_{n-1}}{x_i^l=x_{i+1}^l\ (0\leq i<n)}\]
(subscripts mod~$n$). By relabelling the generators, we see that the group $F^{0/l}(n)$ is isomorphic to $G(n,l)$ if $n$ is odd and is isomorphic to the free product of two copies of $G(n/2,l)$ if $n$ is even. In this context we record the following:
\begin{lemma}\label{lem:F^(0/1)(n)}
Let $G=G(n,l)$ where $n\geq 2, l\geq 1$. Then there is a central extension $\Z\cong \gpres{x_0^l} \hookrightarrow G(n,l) \twoheadrightarrow \underbrace{\Z_l*\cdots *\Z_l}_n$.
\end{lemma}

\begin{proof}
Let $H=\gpres{x_0^l}$, the subgroup of $G$, generated by $x_0^l$. Now, for each $0\leq j<n$, $x_0^lx_jx_0^{-l}x_j^{-1}=x_j^lx_jx_j^{-l}x_j^{-1}=1$, so $x_0^l\in Z(G)$, the centre of $G$. We have $G/H\cong \pres{x_0,\ldots, x_{n-1}}{x_0^l=\cdots=x_{n-1}^l=1}$, which is isomorphic to the free product of $n$ copies of $\Z_l$, and there is an epimorphism $G\rightarrow \Z$ given by sending each $x_i$ to some fixed generator of $\Z$ so $x_i$ (and in particular $x_0$) has infinite order, so $H\cong \Z$.
\end{proof}

\begin{lemma}\label{lem:F^k/lsigns}
For each $k,l\neq 0$ and $n\geq 2$ we have  $F^{k/l}(n)\cong F^{(-k)/(-l)}(n) \cong F^{k/(-l)}(n) \cong F^{(-k)/l}(n)$.
\end{lemma}

\begin{proof}
The isomorphism $F^{k/l}(n)\cong F^{(-k)/(-l)}(n)$ is obtained by replacing each generator by its inverse. We now show that $F^{k/(-l)}(n)\cong F^{(-k)/(-l)}(n)$; the final isomorphism $F^{(-k)/l}(n) \cong F^{k/l}(n)$ is then obtained from this by replacing each generator by its inverse.

The relations of $F^{k/(-l)}(n)$ are $x_i^{-l}x_{i+1}^k=x_{i+2}^{-l}$, which are equivalent to $x_{i+1}^kx_{i+2}^{l}=x_{i}^{l}$. Negating the subscripts and writing $j=-i$ these become $x_{j-1}^kx_{j-2}^{l}=x_{j}^{l}$; adding 2 to the subscripts gives $x_{j+1}^kx_{j}^{l}=x_{j+2}^{l}$. Inverting the relations gives $x_{j}^{-l}x_{j+1}^{-k}=x_{j+2}^{-l}$ which are the relations of $F^{(-k)/(-l)}(n)$, as required.
\end{proof}

Lemmas~\ref{lem:easyFklgroups}--\ref{lem:F^k/lsigns} allow us to assume $k,l\geq 1$.

For our next lemma, recall that a group is \em large \em if it has a finite index subgroup that maps onto the free group of rank~2, that  a group mapping onto a large group is large, and that the free product of two non-trivial finite groups is large unless both groups have order~2 \cite{PrideLargeness}.

\begin{lemma}\label{lem:k/llarge}
Let $n\geq 2$. For each $k,l\geq 1$ let $d=(k,l)$. If $d>1$ then $F^{k/l}(n)$ is large unless $k=n=2$, in which case $F^{k/l}(n)\cong D_\infty$, the infinite dihedral group.
\end{lemma}

\begin{proof}
By killing $x_i^d$ for each $i$ we see that the group $F^{k/l}(n)$ maps onto the free product of $n$ copies of $\Z_d$. Thus $F^{k/l}(n)$ is large if $d>1$ except possibly if $d=2$ and $n=2$, in which case $F^{k/l}(n)\cong \Z_k*\Z_k$ by Lemma~\ref{lem:easyFklgroups}, which is large, unless $k=1$ or $2$. If $k=1$ then $d=1$, a contradiction; if $k=2$ then $F^{k/l}(n)\cong \Z_2*\Z_2=D_\infty$.
\end{proof}

In the notation and terminology of \cite[Chapter~5]{McDermott}, writing $d=(k,l)$, we may express $F^{k/l}(n)=G_n(x_0^lx_1^kx_2^{-l})$ as a composite $G_n(v\circ u$) where $u=x_0^d$, and $v=x_0^{l/d}x_1^{k/d}x_2^{-l/d}$. Since $u$ is a positive word, by \cite[Lemma~5.1.3.4]{McDermott} we then have $G_n(v)=F^{(k/d)/(l/d)}(n)$ embeds in $G_n(v\circ u)=F^{k/l}(n)$. We record this as:

\begin{theorem}[{\cite[Chapter~5]{McDermott}}]\label{thm:F^k/d/l/dembeds}
For each $k,l\geq 1$ let $d=(k,l)$. Then $F^{(k/d)/(l/d)}(n)$ embeds in $F^{k/l}(n)$.
\end{theorem}

In the following corollary, and throughout this paper, by a \em 3-manifold group \em we mean the fundamental group of a (not necessarily closed, compact, or orientable) 3-manifold.

\begin{corollary}\label{cor:embedding}
Let $n\geq 2$, $k,l\geq 1$ and define $d=(k,l)$.
\begin{itemize}
  \item[(a)] Suppose $d>1$. If $F^{(k/d)/(l/d)}(n)$ is not torsion-free then $F^{k/l}(n)$ is an infinite group that is not torsion-free; in particular, if $F^{(k/d)/(l/d)}(n)$ is a finite non-trivial group then $F^{k/l}(n)$ is an infinite group that is not torsion-free.
  \item[(b)] Suppose $F^{(k/d)/(l/d)}(n)$ is not a 3-manifold group; then $F^{k/l}(n)$ is not a 3-manifold group.
  \item[(c)] Suppose $F^{(k/d)/(l/d)}(n)$ is not the fundamental group of an orientable hyperbolic 3-orbifold of finite volume; then $F^{k/l}(n)$ is not the fundamental group of an orientable hyperbolic 3-orbifold of finite volume.
\end{itemize}
\end{corollary}

\begin{proof}
  (a) Since $F^{(k/d)/(l/d)}(n)$ is not torsion-free, it contains a non-trivial element of finite order, which is also an element of $F^{k/l}(n)$. (b) This holds since subgroups of 3-manifold groups are 3-manifold groups \cite[Chapter~8]{Hempel}. (c) If $F^{(k/d)/(l/d)}(n)$ is not the fundamental group of an orientable hyperbolic 3-orbifold of finite volume, then there is no embedding of $F^{(k/d)/(l/d)}(n)$ into $PSL(2,\mathbb{C})$, the group of orientation preserving isometries of hyperbolic 3-space, and hence there is no embedding of $F^{k/l}(n)$ into $PSL(2,\mathbb{C})$, so $F^{k/l}(n)$  is not the fundamental group of an orientable hyperbolic 3-orbifold of finite volume.
\end{proof}

We say that a group $G$ is a \em $q$-generator group \em ($q\geq 1$), or that $G$ is \em $q$-generated, \em if it has a generating set with $q$ generators. Starting with the case $l=1$ we have:

\begin{lemma}\label{lem:Fk2gen}
Let $n\geq 2$, $k\geq 1$. Then $F^{k}(n)$ is 2-generated and can be generated by $x_0$ and $x_1$.
\end{lemma}

\begin{proof}
The relations $x_{i+2}=x_ix_{i+1}^k$ allow each generator $x_{j}$ ($2\leq j< n$) to be written in terms of $x_{j-1}$ and $x_{j-2}$, so only generators $x_0,x_1$ are needed.
\end{proof}

For the general case we have:

\begin{lemma}\label{lem:(k,l)=1isn+1/2gen}
Let $n\geq 3$ be odd, $k,l\geq 1$, $(k,l)=1$. Then $F^{k/l}(n)$ is $(n+1)/2$-generated. In particular, $F^{k/l}(3)$ is 2-generated and can be generated by $x_0$ and $x_1$.
\end{lemma}

\begin{proof}
Since $(k,l)=1$ there exist $\alpha,\beta \in \Z$ such that $\alpha k+\beta l=1$.  The defining relations imply $x_{j+2}^l=x_{j}^lx_{j+1}^k$ and $x_{j+2}^k=x_{j+1}^{-l}x_{j+3}^l$. Thus
\[ x_{j+2}=x_{j+2}^{\alpha k+\beta l}= (x_{j+2}^k)^\alpha (x_{j+2}^l)^\beta = (x_{j+1}^{-l}x_{j+3}^l)^\alpha (x_{j}^lx_{j+1}^k)^\beta\]
and so each generator $x_{j+2}$ can be written in terms of $x_j,x_{j+1},x_{j+3}$. We may therefore eliminate generators $x_{n-1},x_{n-3},\ldots ,x_2$ in turn to leave a presentation with the $(n+1)/2$ generators $x_0,x_1,x_3\ldots ,x_{n-2}$.
\end{proof}

For the case $k=1$ we can decrease the lower bound slightly:

\begin{lemma}\label{lem:F1/latmost(n-1)/2gen}
Let $n\geq 3$, $l\geq 1$. If $n=3$ then $F^{1/l}(n)$ is 2-generated, and if $n\geq 4$ then $F^{1/l}(n)$ is $\lfloor n/2 \rfloor$-generated.
\end{lemma}

\begin{proof}
The relations $x_i^lx_{i+1}=x_{i+2}^l$ can be written $x_{i+1}=x_i^{-l}x_{i+2}^l$. If $n$ is even, this allows all odd numbered generators to be eliminated, leaving an $n/2$-generator presentation. Suppose then that $n$ is odd. Then we can eliminate $x_{2j+1}=x_{2j}^{-l}x_{2(j+1)}^l$ for each $0\leq j\leq (n-3)/2$, leaving a presentation with $(n+1)/2$ generators $x_0,x_2,\ldots ,x_{n-1}$. In doing so, the original relations $x_{n-2}^lx_{n-1}=x_0^l$ and $x_{n-1}^lx_0=x_1^l$ become
\begin{alignat}{1}
(x_{n-3}^{-l}x_{n-1}^l)^lx_{n-1}&=x_0^l,\label{eq:n/2gen1}\\
x_{n-1}^lx_0&=(x_0^{-l}x_2^l)^{l}\label{eq:n/2gen2}.
\end{alignat}
We may substitute the expression for $x_0^l$ given by~(\ref{eq:n/2gen1}) into~(\ref{eq:n/2gen2}) which can then be used to eliminate $x_0$, leaving an $(n-1)/2$-generator presentation.
\end{proof}

The group $F^{k/l}(n)$ has an automorphism $\theta: x_i\mapsto x_{i+1}$ (subscripts mod~$n$), called the \em shift automorphism \em and the corresponding split extension, called the \em shift extension, \em
\[ E^{k/l}(n) = F^{k/l}(n) \rtimes_\theta \pres{t}{t^n}  \]
has a 2-generator, 2-relator presentation
\[ E^{k/l}(n) = \pres{x,t}{t^n, x^ltx^ktx^{-l}t^{-2}}\]
(which is obtained by rewriting $x_0=x$ and $x_i=t^ixt^{-i}$ for $1\leq i<n$).

We now turn to the groups $F^k(n)$. As remarked in \cite[Remark~1]{MaclachlanReid} determining which groups $F^{k}(n)$ (where $n\geq 3$ and odd) are finite is a challenging problem and it is observed that $|F^{2}(3)|=112$, $|F^{3}(3)|=3528$ and that $F^{2}(5)$ is infinite. Using KBMAG \cite{KBMAG} and the {\tt NewmanInfinityCriterion} (\cite{Newman}) command in GAP \cite{GAP} we can prove certain groups $F^{k/l}(n)$ infinite. For example, we have the following result (further infinite groups will be exhibited in Example~\ref{ex:k/l(n;q)}).

\begin{lemma}\label{lem:finiteFk(n)}
Let $n\in\{5,7,9\}$, $3\leq k\leq 12$. Then $F^k(n)$ is infinite.
\end{lemma}

\begin{proof}
If $(n,k)\neq (9,3)$ the group $F^k(n)$ can be proved (automatic and) infinite using KBMAG. The group $F^{3}(9)$ maps onto $H=\pres{x_0,\ldots,x_8}{x_ix_{i+1}^3=x_{i+2}, x_i^{108}\ (0\leq i<9)}$ which can be proved infinite using the {\tt NewmanInfinityCriterion} command applied to the second derived subgroup $H''$ of $H$, with the prime $p=7$.
\end{proof}

\section{Abelianisations}\label{sec:abelianisation}

Knowledge of the order of the abelianisation $|F^{k/l}(n)^\mathrm{ab}|$ will be crucial to our later methods. In Theorem~\ref{thm:Fk/labVn} we obtain a recurrence relation formula for this order. A version of this was asserted in \cite[Lemma, page 238]{SVHigDim} but the formula there is not quite right (for instance, it incorrectly implies that $|F(n)^\mathrm{ab}|$ is even whenever $n$ is odd). While this has no impact on the later arguments in \cite{SVHigDim}, a correct formula is necessary for our arguments, so we include a proof. In Theorem~\ref{thm:Vnsummationformula} we express the order $|F^{k/l}(n)^\mathrm{ab}|$ as a polynomial in $k$ and $l$, and in Corollaries~\ref{cor:Fk/l(n)increasinginn}--\ref{cor:Vnnotdiv(2k)^n} we derive consequences that will be used in later sections.

Define a sequence of natural numbers $V_j^{k/l}$ according to the following recurrence relation
\begin{equation}\label{eq:recrel}
V_1^{k/l}=k,\ V_2^{k/l}=k^2+2l^2,\ V_j^{k/l}=k V_{j-1}^{k/l}+l^2 V_{j-2}^{k/l}\quad (j\geq 3).
\end{equation}

\begin{theorem}\label{thm:Fk/labVn}
Let $k,l\geq 1$, $n\geq 2$. Then
\[ |F^{k/l}(n)^\mathrm{ab}|=\begin{cases}
V_n^{k/l}&\mathrm{if}~n~\mathrm{is~odd};\\
V_n^{k/l}-2l^n&\mathrm{if}~n~\mathrm{is~even}.
\end{cases}\]
\end{theorem}

\begin{proof}
The order $|F^{k/l}(n)^\mathrm{ab}|$ is equal to the resultant $|\mathrm{Res}(f(t),g(t))|$ where $f(t)=l+kt-lt^2$ is the representer polynomial of $F^{k/l}(n)$ and $g(t)=t^n-1$ (see \cite{Johnsonbook}). For each $j\geq 1$ define $u_j=V_j^{k/l}/l^j$ where $V_j^{k/l}$ is as defined at~(\ref{eq:recrel}). Then $f(t)$ is the characteristic polynomial of the recurrence relation defining the sequence $(u_j)$ and has distinct roots $\beta_1,\beta_2$, say. Then the sequence $(u_j)$ has general solution $u_j=c_1\beta_1^j+c_2\beta_2^j$ (see for example Theorems~4.10,4.11 of \cite{NivenZuckermanMontgomery}). Putting $n=1,2$ into these solutions and solving for $c_1,c_2$ gives $c_1=c_2=1$ and hence $u_j=\beta_1^j+\beta_2^j$. Then by \cite[Lemma~2.1]{Odoni}
\[ |\mathrm{Res}(f(t),g(t))|  = |l^n(\beta_1^n-1)(\beta_2^n-1)| = |l^n\left( (-1)^n + 1 -u_n \right)| = |l^n+(-l)^n- V^{k/l}_n|\]
as required.
\end{proof}

This implies, for example, that
\begin{alignat}{1}
|F^{k/l}(3)|=k^3+3kl^2,\label{eq:Fk/l(3)ab}
\end{alignat}
and that $F^{k/l}(n)$ is trivial if and only if $n=2$ and $k=1$.

\begin{theorem}\label{thm:Vnsummationformula}
Let $n,k,l\geq 1$ and let $N=\floor*{n/2}$. Then
$$V_n^{k/l} =\sum_{r=0}^{N}a_{n,r}k^{n-2r}l^{2r}$$
for integers $a_{n,r}\geq 1$ satisfying $a_{n,0}=1$ for $n\geq 1$ and
\( a_{n,r}=a_{n-1,r}+a_{n-2,r-1} \)
for $1\leq r <N$, $n\geq 3$ and $a_{n,N}=n$ if $n$ is odd and $a_{n,N}=2$ if $n$ is even.
\end{theorem}

\begin{proof}
By the definition of $V^{k/l}_n$, the statement is true for $n=1,2$. Suppose that $n\geq 3$ and that the statement is true for all $3\leq j<n$. If $n$ is odd then
	\begin{align*}
	V_n^{k/l}&=kV_{n-1}^{k/l}+l^2V_{n-2}^{k/l}\\
	&=k\left(\sum_{r=0}^{(n-1)/2}a_{n-1,r}k^{n-1-2r}l^{2r}\right)+l^2\left(\sum_{r=0}^{(n-3)/2}a_{n-2,r}k^{n-2-2r}l^{2r}\right)\\
	&=\left( \sum_{r=0}^{(n-1)/2-1}a_{n-1,r}k^{n-2r}l^{2r}+	\sum_{r=0}^{(n-3)/2-1}a_{n-2,r}k^{n-2-2r}l^{2+2r} \right)+nkl^{n-1}\\
	&=\left( \sum_{r=0}^{(n-3)/2}a_{n-1,r}k^{n-2r}l^{2r}+	\sum_{r=1}^{(n-3)/2}a_{n-2,r-1}k^{n-2r}l^{2r} \right)+nkl^{n-1}\\
	&=a_{n-1,0}k^{n}+	\left(\sum_{r=1}^{(n-3)/2}(a_{n-1,r}+a_{n-2,r-1})k^{n-2r}l^{2r}\right) +nkl^{n-1}\\
	&=\sum_{r=0}^{(n-1)/2}a_{n,r}k^{n-2r}l^{2r}
	\end{align*}
	where $a_{n,0}=a_{n-1,0}=1$, $a_{n,(n-1)/2}=n$ and $a_{n,r}=a_{n-1,r}+	a_{n-2,r-1}$ for $1\leq r\leq (n-3)/2$. A similar argument applies when $n$ is even.
\end{proof}

\begin{corollary}\label{cor:Fk/l(n)increasinginn}
Let $k,l\geq 1$. Then for each $n\geq 2$ we have $|F^{k/l}(n+1)^\mathrm{ab}|>|F^{k/l}(n)^\mathrm{ab}|$. Hence if either $n\geq 3$ or ($n=2$ and $k>1$) then the shift automorphism $\theta$ of $F^{k/l}(n)$ has order $n$.
\end{corollary}

\begin{proof}
If $n$ is even then (since, by~(\ref{eq:recrel}), $(V_j)$ is increasing in $j$) we have $|F^{k/l}(n+1)^\mathrm{ab}|=V_{n+1}^{k/l}>V_{n}^{k/l}>V_{n}^{k/l}-2l^{n}= |F^{k/l}(n)^\mathrm{ab}|$. If $n$ is odd then
\[
|F^{k/l}(n+1)^\mathrm{ab}|
= V_{n+1}^{k/l}-2l^{n+1}
= \sum_{r=0}^{(n-1)/2} a_{n+1,r}k^{n+1-2r}l^{2r}
> \sum_{r=0}^{(n-1)/2} a_{n,r}k^{n-2r}l^{2r}=|F^{k/l}(n)^\mathrm{ab}|
\]
since $a_{n+1,0}=a_{n,0}$ and $a_{n+1,r}>a_{n,r}$ for any $r\geq 1$. Now let $n\geq 3$ and suppose that $\theta$ has order $m|n$. If $m=1$ then $F^{k/l}(n)\cong \Z_k$, which contradicts $|F^{k/l}(n)^\mathrm{ab}|\geq k^2+2l^2$; if $2\leq m<n$ then $F^{k/l}(n)\cong F^{k/l}(m)$ and so $|F^{k/l}(n)^\mathrm{ab}|=|F^{k/l}(m)^\mathrm{ab}|$, a contradiction. Finally, if $n=2$, $k>1$, and $m=1$ then $F^{k/l}(n)\cong \Z_k$, which contradicts $|F^{k/l}(2)^\mathrm{ab}|=k^2$.
\end{proof}

\begin{corollary}\label{cor:Vntilde}
Suppose $n\geq 3$ is odd, $k,l\geq 1$, $(k,l)=1$ where $k$ is even. Then the 2-adic orders $v_2(k)$ and $v_2(V_n^{k/l})$ are equal.
\end{corollary}

\begin{proof}
Let $k=2^mq$ where $q$ is odd and $m\geq 1$. We claim that for each $n\geq 1$ there exists some odd $\tilde{V}_n^{k/l}$ such that $V_n^{k/l} = 2^m\tilde{V}^{k/l}_n$ if $n$ is odd and $V_n^{k/l} = 2\tilde{V}^{k/l}_n$ if $n$ is even.

 As at~(\ref{eq:recrel}) we have $V_1^{k/l}=k=2^mq=$ $2^m\tilde{V}_1^{k/l}$ where $\tilde{V}_1^{k/l}=q$ is odd, 	 $V_2^{k/l}=k^2+2l^2= $ $2\tilde{V}_2^{k/l}$ where $\tilde{V}_2^{k/l}=k^2/2+l^2$, which is odd.
Suppose that $n\geq 3$ and that the claim holds  for all $j<n$. If $n$ is odd, then
 \[ V_n^{k/l}=kV_{n-1}^{k/l}+l^2V_{n-2}^{k/l}=k(2\tilde{V}_{n-1}^{k/l})+l^2(2^m\tilde{V}_{n-2}^{k/l})= 2^m \tilde{V}^{k/l}_n \]
where $\tilde{V}^{k/l}_n=2q\tilde{V}_{n-1}^{k/l}+l^2\tilde{V}_{n-2}^{k/l}$ is odd.
Similarly if $n$ is even, then
\[ V_n^{k/l}=kV_{n-1}^{k/l}+l^2V_{n-2}^{k/l}=k(2^m\tilde{V}_{n-1}^{k/l})+l^2(2\tilde{V}_{n-2}^{k/l})= 2 \tilde{V}^{k/l}_n \]
where $\tilde{V}^{k/l}_n=2^{m-1}k\tilde{V}_{n-1}^{k/l}+l^2\tilde{V}_{n-2}^{k/l}$ is odd.
\end{proof}

\begin{corollary}\label{cor:Vneven}
Let  $n,k,l\geq 1$. Then $V_n^{k/l}$ is even if and only if either $k$ is even or ($l$ is odd and $n\equiv 0$~mod~$3$).
\end{corollary}

\begin{proof}
Note that $F^{k/l}(n)$ maps onto $\Z_k$ (by sending each $x_i$ to some fixed generator of $\Z_k$) so if $k$ is even then $|F^{k/l}(n)^\mathrm{ab}|$ is even, so assume $k$ is odd. If $l$ is even then by~(\ref{eq:recrel}), $V_n^{k/l}\equiv V_{n-1}^{k/l}\bmod 2$ for all $n\geq 2$ and $V_1^{k/l}=k$ is odd, and so $V_n^{k/l}$ is odd for all $n\geq 1$. If $l$ is odd then $V_1^{k/l}$ and $V_2^{k/l}$ are odd, and $V_n^{k/l}\equiv V_{n-1}^{k/l}+V_{n-2}^{k/l}\bmod 2$ for all $n\geq 3$ and it follows that $V_n^{k/l}$ is even if and only if $n\equiv 0\bmod 3$.
\end{proof}

\begin{corollary}\label{cor:Fk/ldivides}
  Let $n\geq 3$ be odd, $k,l\geq 1$, and suppose that $(k,l) = 1$.
  \begin{itemize}
    \item[(a)] If $|F^{k/l}(n)^\mathrm{ab}|$ divides $(2l)^n$ then $k=l=1$ and $n=3$;
    \item[(b)] if $|F^{k/l}(n)^\mathrm{ab}|$ divides $(2k)^n$ then $|F^{k/l}(n)^\mathrm{ab}|$ is even and $k$ is odd.
  \end{itemize}
\end{corollary}

\begin{proof}
(a) First we claim that for each $n\geq 1$ we have $(V_n^{k/l},l)=1$ (which, by definition of $V_n^{k/l}$, is true for $n=1,2$). Suppose this statement is true for $j-1$ where $j\geq 3$. Then
\[(V_j^{k/l},l)=(kV_{j-1}^{k/l}+l^2V_{j-2}^{k/l},l)=(kV_{j-1}^{k/l},l)=1\]
so by induction, the statement is true for all $n\geq 1$. Therefore by Theorem~\ref{thm:Fk/labVn} we have $(|F^{k/l}(n)^{\mathrm{ab}}|,l)=(V_n^{k/l},l)=1$ for all $n\geq 1$.

Suppose that $|F^{k/l}(n)^{\mathrm{ab}}|$  divides $(2l)^n$. Then $|F^{k/l}(n)^{\mathrm{ab}}|$ divides $2^n$ since $(|F^{k/l}(n)^{\mathrm{ab}}|,l)=1$.
By Theorem \ref{thm:Fk/labVn} we have $|F^{k/l}(n)^{\mathrm{ab}}|=V_n^{k/l}$ and by Theorem~\ref{thm:Vnsummationformula}
\[V_n^{k/l}=k^{n}+	\left(\sum_{r=1}^{N-1}a_{n,r}k^{n-2r}l^{2r}\right) +nkl^{n-1}\]
where $N=\floor*{\frac{n}{2}}$ and each $a_{n,r}\geq 1$ is an integer so in particular, $k^n<2^n$  and $nkl^{n-1} < 2^n$ and hence $k=l=1$. Then $V_j^{k/l}$ is a Lucas number, which therefore divides $2^n$ and since the only powers of 2 that appear in the Lucas sequence are 1,2,4 (see, for example, \cite{BravoLuca}) we have $n=3$.

(b) Suppose that $|F^{k/l}(n)^{\mathrm{ab}}|$ divides $(2k)^n$. If $|F^{k/l}(n)^{\mathrm{ab}}|$ is odd, then it divides $k^n$ which is impossible since $V_n^{k/l}>k^n$ by Theorem~\ref{thm:Vnsummationformula}. Therefore $|F^{k/l}(n)^{\mathrm{ab}}|$ is even.  Suppose for contradiction that $k$ is even, say $k=2^mq$ where $q$ is odd and $m\geq 1$. Then by Theorem~\ref{thm:Fk/labVn} and Corollary~\ref{cor:Vntilde} we have $|F^{k/l}(n)^{\mathrm{ab}}|/2^m$ is odd, and so $|F^{k/l}(n)^{\mathrm{ab}}|/2^m$ divides $q^n$. But by Theorem~\ref{thm:Vnsummationformula}
\[|F^{k/l}(n)^{\mathrm{ab}}|/2^m= a_{n,0}q^{n}2^{m(n-1)}+	\sum_{r=1}^{N}a_{n,r}k^{n-2r}l^{2r}2^{-m} >q^n\]
since $N\geq 1$, a contradiction. Therefore $k$ is odd.
\end{proof}

\begin{corollary}\label{cor:Vnnotdiv(2k)^n}
Suppose $n=pk>7$ is odd, where $p\geq 1$, $k\geq 3$, $(p,k)=1$, and $l\geq 1$. If $(k,l)=1$ then $V_n^{k/l}$ does not divide $(2k)^n$.
\end{corollary}

\begin{proof}
Suppose for contradiction that $V_n^{k/l}$ divides $(2k)^n$. By Theorem~\ref{thm:Vnsummationformula} we have
\[ V_n^{k/l} = k^2\left( \left( k\sum_{r=0}^{(n-3)/2} a_{n,r} k^{n-2r-3} l^{2r} \right) + pl^{n-1}\right)\]
so $V_n^{k/l}\equiv 0$~mod~$k^2$ and $(V_n^{k/l}/k^2,k)=1$ and so $V_n^{k/l}/k^2$ divides $2^n$. But
\begin{alignat*}{1}
\frac{V_n^{k/l}}{k^2} &= \left( k \sum_{r=0}^{(n-3)/2} a_{n,r} k^{n-2r-3}l^{2r}\right) + pl^{n-1}\\
&> k^{n-2}+k^{n-4}\geq k^{n-4}(3^2 +1)=10k^{n-4}>2^n
\end{alignat*}
since $n\geq 7$ and $k\geq 3$, a contradiction.
\end{proof}

\section{Hyperbolic 3-orbifolds}\label{sec:hyporbifolds}

In this section we prove the following.

\begin{theorem}\label{thm:hyper}
Let $n,k,l\geq 1$, where $n$ is odd. Then $F^{k/l}(n)$ is not the fundamental group of an orientable hyperbolic 3-orbifold (in particular, 3-manifold) of finite volume.
\end{theorem}

\begin{corollary}\label{cor:hyperkodd}
Let $n,k,l\geq 1$, where $n$ and $k$ are odd. Then $F^{k/l}(n)$ is not the fundamental group of a hyperbolic 3-orbifold (in particular, 3-manifold) of finite volume.
\end{corollary}

Note that we do not assume $(k,l)=1$ in the hypotheses of Theorem~\ref{thm:hyper} and Corollary~\ref{cor:hyperkodd}. Our method of proof follows that introduced in \cite{Maclachlan} (for Fibonacci groups $F(n)$), and developed further in \cite{BardakovVesnin,CRStopprop,SVFibGpsoddnumbergens}. That is, supposing that $F^{k/l}(n)$ is the fundamental group of an orientable hyperbolic 3-orbifold of finite volume, then so is its shift extension $E^{k/l}(n)=\pres{x,t}{t^n, x^ltx^ktx^{-l}t^{-2}}$, which is therefore isomorphic to a subgroup of $PSL(2,\mathbb{C})$. We show that a putative embedding in $PSL(2,\mathbb{C})$ would imply restrictions on the order of the abelianisation $F^{k/l}(n)^\mathrm{ab}$ and then use the results of Section~\ref{sec:abelianisation} to show that these restrictions cannot occur.

\begin{proof}[{Proof of Theorem~\ref{thm:hyper}}]
We prove the theorem in the case $(k,l)=1$; the case $(k,l)>1$ then follows from Corollary~\ref{cor:embedding}. If $n=1$ then $F^{k/l}(n)\cong \Z_k$, so assume $n\geq 3$. Suppose for contradiction that $F^{k/l}(n)$ is the fundamental group of an orientable hyperbolic 3-orbifold of finite volume. By Corollary~\ref{cor:Fk/l(n)increasinginn} the shift automorphism $\theta$ of $F^{k/l}(n)$ has order $n$ so, as explained in the proof of \cite[Theorem~3.1]{Maclachlan}, it follows from the Mostow Rigidity Theorem that the shift extension $E=\pres{x,t}{t^n, x^ltx^ktx^{-l}t^{-2}}$ of $F^{k/l}(n)$ is isomorphic to a subgroup of $PSL(2,\mathbb{C})$.

Therefore there exists a subgroup $\tilde{E}$  of $SL(2,\mathbb{C})$, which is the pre-image of $E$ with respect to the canonical projection. Suppose that, for the generator $x$ of $E$, the corresponding element in $\tilde{E}$ is the matrix $\tilde{x}=\begin{bmatrix}
	a       & b  \\
	c       & d
	\end{bmatrix}\in SL(2,\mathbb{C})$, where (as in the proof of \cite[Theorem~3.1]{CRStopprop}) $bc \neq 0$, since $E$ has finite covolume.
For the generator $t\in E$, the corresponding element in $\tilde{E}$ is the matrix $\tilde{t}=\begin{bmatrix}
	\zeta       & 0  \\
	0       & \zeta^{-1}
	\end{bmatrix}\in SL(2,\mathbb{C})$, where $\zeta$ is a primitive root of unity in $\mathbb{C}$ of order $2n$. Then the relation
\begin{align*}\label{EQ1}
	tx^kt^{-1}=x^{-l}t^2x^lt^{-2}
\end{align*}
induces the relation
	\begin{equation}
	\begin{bmatrix}
	\epsilon       & 0  \\
	0       & \epsilon
	\end{bmatrix}\begin{bmatrix}
	\zeta       & 0  \\
	0       & \zeta^{-1}
	\end{bmatrix}
	\begin{bmatrix}
	a       & b  \\
	c       & d
	\end{bmatrix}^k
	\begin{bmatrix}
	\zeta      & 0  \\
	0       & \zeta^{-1}
	\end{bmatrix}^{-1}
	=\begin{bmatrix}
	d       & -b  \\
	-c       & a
	\end{bmatrix}^l\begin{bmatrix}
	\zeta       & 0  \\
	0       & \zeta^{-1}
	\end{bmatrix}^2	
	\begin{bmatrix}
	a       & b  \\
	c       & d
	\end{bmatrix}^l\begin{bmatrix}
	\zeta       & 0  \\
	0       & \zeta^{-1}
	\end{bmatrix}^{-2}	\label{eq:eqq}
	\end{equation}
	where $\epsilon=\pm 1$. It was observed in \cite[page~962]{SVFibGpsoddnumbergens} that in $SL(2,\mathbb{C})$
\begin{align*}
	\begin{bmatrix}
	 a       & b  \\
	c        &  d
	\end{bmatrix}^j&=\begin{bmatrix}
	 S_j     &  bR_j \\
	cR_j        & T_j
	\end{bmatrix}
	\end{align*}
	where $S_{j+1}=aS_j+bcR_j$, $T_{j+1}=dT_j+bcR_j$, and $R_{j+1}=S_j+dR_j$, with $S_1=a$, $T_1=d$ and $R_1=1$. Note that the determinant $S_jT_j-bcR_j^2=1$. Applying this formula to our case, the left hand side of~(\ref{eq:eqq}) is
	\begin{align}\label{eq:EQ2}
	\begin{bmatrix}
	\epsilon       & 0  \\
	0       & \epsilon
	\end{bmatrix}\begin{bmatrix}
	\zeta       & 0  \\
	0       & \zeta^{-1}
	\end{bmatrix}
	\begin{bmatrix}
	S_k       & bR_k \\
	c R_k       & T_k
	\end{bmatrix}
	\begin{bmatrix}
	\zeta      & 0  \\
	0       & \zeta^{-1}
	\end{bmatrix}^{-1}
	&=\begin{bmatrix}
	\epsilon S_k      & \epsilon \zeta^2 b R_k \\
	\epsilon \zeta^{-2} c R_k       & \epsilon T_k
	\end{bmatrix}.
	\end{align}
Similarly,  the right hand side of~(\ref{eq:eqq}) is
\begin{align}
  \begin{bmatrix}
	T_l       & -bR_l  \\
	-c R_l      & S_l
	\end{bmatrix}
    \begin{bmatrix}
	\zeta^2       & 0  \\
	0       & \zeta^{-2}
	\end{bmatrix}	
	\begin{bmatrix}
	S_l       & bR_l  \\
	c R_l     & T_l
	\end{bmatrix}
    \begin{bmatrix}
	\zeta^{-2}       & 0  \\
	0       & \zeta^2
	\end{bmatrix}
	=\begin{bmatrix}
	T_lS_l-\zeta^{-4}bcR_l^2       & (\zeta^4-1)bT_lR_l  \\
	(\zeta^{-4}-1)cS_lR_l       & T_lS_l-\zeta^{4}bcR_l^2
	\end{bmatrix}.	\label{eq:EQ2again}
\end{align}
Therefore, since $T_lS_l-bcR_l^2=1$, the equations~(\ref{eq:eqq}),(\ref{eq:EQ2}),(\ref{eq:EQ2again}) give
\begin{align}
	\begin{bmatrix}
	\epsilon S_k      & \epsilon b\zeta^2 R_k \\
	\epsilon \zeta^{-2} cR_k       & \epsilon T_k
	\end{bmatrix}&=\begin{bmatrix}
	T_lS_l-\zeta^{-4}bcR_l^2       & (\zeta^4-1)bT_lR_l  \\
	(\zeta^{-4}-1)cS_lR_l       & T_lS_l-\zeta^{4}bcR_l^2
	\end{bmatrix}\nonumber\\
	&=\begin{bmatrix}\label{neq1}
	1+(1-\zeta^{-4})bcR_l^2       & (\zeta^4-1)bT_lR_l  \\
	(\zeta^{-4}-1)cS_lR_l       & 1+(1-\zeta^{4})bcR_l^2
	\end{bmatrix}.
	\end{align}

Comparing the terms on both sides of  (\ref{neq1}) gives
\begin{align*}
\epsilon b \zeta^2  R_k&=(\zeta^4-1)bT_lR_l,\\
\epsilon c \zeta^{-2}  R_k  &=(\zeta^{-4}-1)cS_lR_l.
\end{align*}
Since $bc\neq 0$, we get
\begin{align}
\epsilon R_k&=(\zeta^2 -\zeta^{-2})R_lT_l,\label{eq:Eqk1}\\
\epsilon R_k  &=-(\zeta^2 -\zeta^{-2})R_lS_l,\nonumber
\end{align}
and, since $\zeta$ is a primitive root of unity of order $2n$ with $n$ odd $\zeta^2 -\zeta^{-2}\neq 0$, so
\begin{align}\label{eq:EQ3}
R_l(T_l+S_l)=0.
\end{align}

Suppose $R_l=0$. Then $\epsilon T_k=\epsilon S_k=1$ by (\ref{neq1}) and  $R_k=0$ by (\ref{eq:Eqk1}). Hence $\tilde{x}^k=\begin{bmatrix}
S_k      & b R_k \\
 cR_k       &  T_k
\end{bmatrix}=\pm I$, so $x_i^{2k}=1$ for all generators $x_i$ of $F^{k/l}(n)$. Thus the order $|F^{k/l}(n)^{\mathrm{ab}}|$ divides $(2k)^n$ and then Corollary~\ref{cor:Fk/ldivides}
implies that $|F^{k/l}(n)^\mathrm{ab}|$ is even and $k$ is odd. Then by Corollary~\ref{cor:Vneven} $l$ is odd and $n\equiv 0$~mod~$3$. The map from $F^{k/l}(n)$ to $F^{k/l}(3)$ (and hence from $F^{k/l}(n)^{\mathrm{ab}}$ to $F^{k/l}(3)^{\mathrm{ab}}$) sending $x_i$ to $x_{i\bmod 3}$ is a surjective homomorphism. Hence $x_0$,$x_1$ and $x_2$ each has order dividing $2k$ in $(F^{k/l}(3))^{\mathrm{ab}}$ and so $|F^{k/l}(3)|^{\mathrm{ab}}$ divides $(2k)^3$. But $|F^{k/l}(3)|^{\mathrm{ab}}=k^3+3kl^2$ which divides $(2k)^3$, and so $k^2+3l^2$ divides $8k^2$. Therefore there exists a natural number $m$ such that
\begin{alignat}{1}
l^2=\dfrac{(8-m)k^2}{3m}\label{eq:lsquared}
\end{alignat}
and so $m\in \lbrace 1,2,\ldots ,7 \rbrace$. When  $m=1,2,3,4,5,6, 7$, equation (\ref{eq:lsquared}) implies $k=\sqrt{3}l/\sqrt{7},l,3l/\sqrt{5}, \sqrt{3}l,$ $\sqrt{5}l,3l, \sqrt{21}l$ respectively. Hence $m=2$ or $6$, and so either $k=l=1$ or $k=3$ and $l=1$. If $k=l=1$ then the result is given in \cite[Theorem~3.1]{Maclachlan} so assume $k=3,l=1$, and therefore $F^{3/1}(n)$ divides $(2k)^n=6^n$. If $9|n$ then $|F^{3/1}(9)^\mathrm{ab}|=2^2\cdot 3^3\cdot 433$ divides $|F^{3/1}(n)^\mathrm{ab}|$, a contradiction. Thus $n=3p$ for some $p$ where $(p,3)=1$. If $n\geq 9$ then the result follows from Corollary~\ref{cor:Vnnotdiv(2k)^n} and if $n=3$ then a computation in GAP shows that $F^{3/1}(3)$ is a finite group of order $3528$ which cannot occur as subgroup of $PSL(2,\mathbb{C})$ (see, for example, \cite[pages~152--154]{Lyndon}).

Now suppose that $R_l\neq  0$. Then  $T_l=-S_l$ by (\ref{eq:EQ3}) so $\tilde{x}^l$ is traceless, and so $\tilde{x}^{2l}=-I$, where $I$ is the identity element of $SL(2,\mathbb{C})$. Therefore $x^{2l}=I$ and so $x_i^{2l}=1$ for all generators $x_i$ of $F^{k/l}(n)$. Thus the order $|F^{k/l}(n)^{\mathrm{ab}}|$ divides $(2l)^n$. Hence by Corollary~\ref{cor:Fk/ldivides} $k=l=1$, and $n=3$, in which case $F^{k/l}(n)=F(3)\cong Q_8$, which is not the fundamental group of a hyperbolic 3-orbifold.
\end{proof}

\begin{proof}[{Proof of Corollary~\ref{cor:hyperkodd}}]
Since $k$ is odd, the defining relators of $F^{k/l}(n)$ imply that for each $0\leq i<n$ generator $x_{i+1}=(x_{i+1}^{-(k-1)/2})^2 x_i^{-l}x_{i+2}^l$, which is a product of an even number of generators. Hence if $G$ is the fundamental group of a hyperbolic 3-orbifold of finite volume, then that orbifold must be orientable, which is not possible by Theorem~\ref{thm:hyper}.
\end{proof}

\section{Torsion and asphericity}\label{sec:torsion}

In this section we fix $w(n,k)=x_0^kx_1^k\ldots x_{n-1}^k \in F^{k/l}(n)$. Our starting point is the following result of Bardakov and Vesnin \cite{BardakovVesnin}, who note that in the case $k=l=1$ the words $w(n,1)\in F(n)$ were first considered by Johnson \cite{Johnsonbook}.

\begin{theorem}[{\cite[Proposition~3.1]{BardakovVesnin}}]\label{thm:Worder2inF(n)}
Suppose $n\geq 9$ is odd. Then $w(n,1)$ is an element of order~2 in (the infinite group) $F(n)$.
\end{theorem}

We have the following corollary concerning the asphericity of the relative presentation of the shift extension of $F^{k/l}(n)$, where the terms \em relative presentation \em and  \em aspherical \em are as defined in \cite{BEW}.

\begin{corollary}[{compare \cite[Example~4.3(a)]{BEW}}]\label{cor:F(n)extnotasp}
Suppose $n\geq 9$ is odd. Then the relative presentation $\pres{G,x}{xtxtx^{-1}t^{-2}}$ (where $G=\pres{t}{t^n}$) is not aspherical.
\end{corollary}

(If $n\in \{3,5,7\}$ then $w(n,1)=1$ in the finite group $F(n)$.) Theorem~\ref{thm:Worder2inF(n)} is significant because it gives examples of infinite cyclically presented groups with torsion. Indeed, in many studies (for example \cite{GilbertHowie,Prishchepov,BardakovVesnin,CRS05,EdjvetWilliams,BogleyParker,Spaggiari}) cyclically presented groups are proved infinite by showing that they are non-trivial and that the relative presentations of their shift extensions are aspherical, and deducing (by \cite[Lemma~3.1]{GilbertHowie}, \cite[Theorem~4.1(a)]{Bogleyshift}) that the cyclic presentation is topologically aspherical, and hence that the cyclically presented group is torsion-free.

In this section we obtain similar results to Theorem~\ref{thm:Worder2inF(n)} and Corollary~\ref{cor:F(n)extnotasp} for groups $F^{k/l}(n)$ under certain conditions on $k,l$. Theorem~\ref{thm:order2inFk/l}(a) generalizes \cite[Exercise~12, page 84]{Johnsonbook} from $F(n)$ to $F^{k/l}(n)$; part (b) generalizes the first part of the proof of \cite[Proposition~3.1]{BardakovVesnin} (see also \cite[Exercise~2, page 83]{Johnsonbook}) from the groups $F(n)$ to the groups $F^{k/l}(n)$. We use the notation $[a,b]=a^{-1}b^{-1}ab$.

\begin{theorem}\label{thm:order2inFk/l}
Let $n\geq 3$ be odd, $k,l\geq 1$ and let $w(n,k)=x_0^kx_1^k\ldots x_{n-1}^k \in F^{k/l}(n)$. Then
\begin{itemize}
  \item[(a)] $w(n,k)^2=1$;
  \item[(b)] $w(n,k)=[x_0^l,x_{n-1}^l]$.
\end{itemize}
\end{theorem}

\begin{proof}
(a)
We have
\begin{alignat*}{1}
w(n,k)^2
&= (x_0^kx_1^k)(x_2^kx_3^k) \ldots (x_{n-1}^kx_0^k) \ldots (x_{n-2}^kx_{n-1}^k)\\
&= (x_0^{-(l-k)} x_2^k x_2^{l-k}) (x_2^{-(l-k)} x_4^k x_4^{l-k}) \ldots (x_{n-1}^{-(l-k)} x_1^k x_1^{l-k}) \ldots (x_{n-2}^{-(l-k)} x_0^k x_0^{l-k})\\
&= x_0^{-l} \left( x_0^k x_2^k x_4^k\ldots x_{n-2}^k \right) x_0^l\\
&= x_0^{-l}  \left( (x_{n-1}^{-l}x_{1}^l) (x_{1}^{-l}x_{3}^l) (x_{3}^{-l}x_{5}^l) \ldots (x_{n-3}^{-l}x_{n-1}^l) \right) x_0^l\\
&= 1.
\end{alignat*}

(b) We have
\begin{alignat*}{1}
w(n,k)
&=x_{0}^kx_{1}^kx_{2}^kx_{3}^kx_{4}^kx_{5}^k\ldots x_{n-1}^k\\
&=x_{0}^{-l}x_{0}^k (x_0^l x_{1}^k) x_{2}^kx_{3}^kx_{4}^kx_{5}^k\ldots x_{n-1}^k\\
&=x_{0}^{-l}x_{0}^k (x_{2}^l) x_{2}^kx_{3}^kx_{4}^kx_{5}^k\ldots x_{n-1}^k\\
&=x_{0}^{-l}x_{0}^k x_{2}^k (x_{2}^l x_{3}^k) x_{4}^kx_{5}^k\ldots x_{n-1}^k\\
&= \cdots \\
&= x_{0}^{-l}x_{0}^k x_{2}^k x_{4}^k x_{6}^k \ldots x_{n-1}^k x_{n-1}^l\\
&= x_{0}^{-l}(x_{n-1}^{-l}x_{1}^l) (x_{1}^{-l}x_{3}^l) (x_{3}^{-l}x_{5}^l) (x_{5}^{-l}x_{7}^l) \ldots (x_{n-2}^{-l}x_{0}^l) x_{n-1}^l\\
&= x_{0}^{-l} x_{n-1}^{-l}x_{0}^l x_{n-1}^l.
\end{alignat*}
\end{proof}

As we now show, in many cases (for odd $n$) we have $w(n,k)\neq 1$, and so $w(n,k)$ is an element of order 2. Corollary~\ref{cor:Fkabelian} generalizes the second part of the proof of \cite[Proposition~3.1]{BardakovVesnin} from the groups $F(n)$ to the groups $F^k(n)$, showing that for odd $n$ the group $F^{k}(n)$ is not torsion-free. This is in contrast to the case where $n$ is even where, if either $k=1$ and $n\geq 8$ or $k\geq 2$ and $n\geq 6$ the group $F^k(n)$ (being the fundamental group of a hyperbolic manifold \cite[Theorem~C]{HKM},\cite[Theorem~3]{MaclachlanReid}) is torsion-free.

\begin{corollary}\label{cor:Fkabelian}
Let $n\geq 3$ be odd, $k\geq 1$, $G=F^{k}(n)$ and let $w(n,k)=x_0^kx_1^k\ldots x_{n-1}^k \in G$. Then the normal closure of $w(n,k)$ in $G$ is equal to the derived subgroup of $G$. Thus $w(n,k)=1$ if and only if $G$ is abelian. In particular, $G$ is not torsion-free.
\end{corollary}

\begin{proof}
By Theorem~\ref{thm:order2inFk/l} $w(n,k)=[x_0,x_1]$, and by Lemma~\ref{lem:Fk2gen} $G$ is generated by $x_0,x_1$ so the derived subgroup $G'=\npres{w}^G$. For the `in particular', note that if $G$ is infinite, then since $G^\mathrm{ab}$ is finite, $w$ is of order 2, so $G$ is not torsion-free, and if $G$ is finite then it is not torsion-free, since it is non-trivial.
\end{proof}

For the case $n=3$ we have the following:

\begin{corollary}\label{cor:Fk/l(3)abelian}
Let $k,l\geq 1$, $(k,l)=1$, $G=F^{k/l}(3)$ and let $w(3,k)=x_0^kx_1^kx_{2}^k \in G$. Then the normal closure of $w(3,k)$ in $G$ is equal to the derived subgroup of $G$. Thus $w(3,k)=1$ if and only if $G$ is abelian. In particular, $G$ is not torsion-free.
\end{corollary}

\begin{proof}
Let $w=w(3,k)$, $N=\npres{w}^G$. By Theorem~\ref{thm:order2inFk/l} $w=[x_0^l,x_1^l]\in G'$, so $N$ is a subgroup of $G'$. We shall show that $G/N$ is abelian, and so $G'$ is a subgroup of $N$, and hence $N=G'$. The `in particular' will follow as in the proof of Corollary~\ref{cor:Fkabelian}.

In $G/N$ we have $x_0^kx_2^kx_1^k=(x_2^{-l}x_1^l)(x_1^{-l}x_0^l)(x_0^{-l}x_2^l)=1$ and $x_0^kx_1^kx_2^k=w(3,k)=1$ and hence $x_1^kx_0^k=x_2^{-k}=x_0^kx_1^k$, $x_2^kx_1^k=x_0^{-k}=x_1^kx_2^k$, $x_0^kx_2^k=x_1^{-k}=x_2^kx_0^k$. That is, $[x_j^k,x_{j+1}^k]=1$ for each $0\leq j\leq 2$. Moreover, by Theorem~\ref{thm:order2inFk/l}(b),  for each $0\leq j<2$ we have $1=\theta^{j+1}(w)=\theta^{j+1}([x_0^l,x_{n-1}^l])=[x_{j+1}^l,x_{j}^l]$.

The relations $x_j^lx_{j+1}^k=x_{j+2}^l$ and $[x_{j}^k,x_{j+2}^k]=1$ in $G/N$ imply
\[x_j^kx_{j+2}^l=x_j^kx_j^lx_{j+1}^k= x_j^lx_j^kx_{j+1}^k=x_j^lx_{j+1}^kx_{j}^k= x_{j+2}^lx_j^k.\]
Hence
\[ x_{j+2}^kx_j^l = x_{j+2}^kx_{j+1}^lx_{j+2}^k = x_{j+1}^l x_{j+2}^k x_{j+2}^k = x_{j}^l x_{j+2}^k.\]
Since $(k,l)=1$ there exist $\alpha,\beta\in \Z$ such that $\alpha k+\beta l=1$. Then for each $j$ we have
\begin{alignat*}{1}
&\ x_jx_{j+2}
=x_j^{\alpha k +\beta l}x_{j+2}^{\alpha k +\beta l}
=x_j^{\alpha k} x_j^{\beta l} x_{j+2}^{\beta l}x_{j+2}^{\alpha k}
=x_j^{\alpha k} x_{j+2}^{\beta l} x_j^{\beta l} x_{j+2}^{\alpha k}
=x_{j+2}^{\beta l} x_j^{\alpha k} x_j^{\beta l} x_{j+2}^{\alpha k}\\
&\quad=x_{j+2}^{\beta l} x_j^{\beta l} x_j^{\alpha k}  x_{j+2}^{\alpha k}
=x_{j+2}^{\beta l} x_j^{\beta l} x_{j+2}^{\alpha k} x_j^{\alpha k}
=x_{j+2}^{\beta l} x_{j+2}^{\alpha k} x_j^{\beta l}  x_j^{\alpha k}
=x_{j+2}^{\alpha k+\beta l} x_j^{\alpha k+ \beta l}
=x_{j+2}x_j.
\end{alignat*}
Hence $G/N$ is abelian.
\end{proof}

\begin{example}\label{ex:F1/2(5)}
\em Let $G=F^{1/2}(5)$. Using GAP we see that $G/\npres{w(5,1)}^G\cong \Z_{101}=G^{\mathrm{ab}}$. Moreover, a computation using KBMAG shows that $G$ is infinite, so non-abelian, and thus $w(5,1)\neq 1$ in $F^{1/2}(5)$, which is therefore not torsion-free.\em
\end{example}

Thus Corollaries~\ref{cor:Fkabelian},\ref{cor:Fk/l(3)abelian} and Example~\ref{ex:F1/2(5)} give cases where $w(n,k)=1$ is equivalent to $F^{k/l}(n)$ being abelian. We expect that in most cases $F^{k/l}(n)$ ($n$ odd) is not abelian, and thus $w(n,k)\neq 1$. However, in some cases $F^{k/l}(n)$ is abelian. The cases we know of are $F(5)\cong \Z_{11}$, $F(7)\cong \Z_{29}$, $F^{1/2}(3)\cong \Z_{13}$ and $F^{2/3}(3)\cong \Z_{62}$. It would be interesting to know if there are any further cases. We know of the following finite non-abelian groups $F^{k/l}(n)$: $F(3)\cong Q_8$;  $F^{1/3}(3)$ of order 3584;  $F^{1/4}(3)$ of order 392;  $F^{2}(3)$ of order 112;  $F^{3}(3)$ of order 3528;  $F^{3/2}(3)$ of order 504. In each of these cases $n=3$ so $w(3,k)\neq 1$ by Corollary~\ref{cor:Fk/l(3)abelian}.

We now turn to the question of asphericity.

\begin{corollary}\label{cor:Fkl(n)extnotasp}
Suppose $n\geq 3$ is odd and let $k,l\geq 1$. If $F^{k/l}(n)$ is finite or $w(n,k)\neq 1$ in $F^{k/l}(n)$ then the relative presentation $\mathcal{P}=\pres{G,x}{x^ltx^ktx^{-l}t^{-2}}$ (where $G=\pres{t}{t^n}$) is not aspherical.
\end{corollary}

\begin{proof}
The group $G(\mathcal{P})$ defined by $\mathcal{P}$ is isomorphic to the shift extension of $F^{k/l}(n)$ so $F^{k/l}(n)$ is isomorphic to a subgroup of $G(\mathcal{P})$. By \cite[Section~3]{Huebschmann} (due to Serre), if the relative presentation $\mathcal{P}$  is aspherical then every finite subgroup of $G(\mathcal{P})$ is conjugate to a subgroup of $G$ (see also \cite[Theorem~2.4(c)]{BEW}). If $F^{k/l}(n)$ is finite and conjugate to a subgroup of $G$ then $F^{k/l}(n)$ is abelian of order at most $n$; but by Theorem~\ref{thm:Fk/labVn} $|F^{k/l}(n)^\mathrm{ab}|\geq |F(n)^\mathrm{ab}|>n$ for all $n\geq 3$, a contradiction. If $w(n,k)\neq 1$ then
 it generates a cyclic subgroup of $F^{k/l}(n)$ of order 2 which, since $n$ is odd, is not conjugate to a subgroup of $G$. Thus $\mathcal{P}$ is not aspherical.
\end{proof}

\begin{corollary}\label{cor:Fk(n)extnotasp}
Suppose $n\geq 3$ is odd and $k\geq 1$. Then the relative presentation $\mathcal{P}=\pres{G,x}{xtx^ktx^{-1}t^{-2}}$ (where $G=\pres{t}{t^n}$) is not aspherical.
\end{corollary}

\begin{proof}
By Corollary~\ref{cor:Fkl(n)extnotasp} we may assume $w(n,k)=1$, so $F^{k}(n)$ is abelian, by Corollary~\ref{cor:Fkabelian}. But then $F^{k/l}(n)$ is finite, so the result follows from Corollary~\ref{cor:Fkl(n)extnotasp}.
\end{proof}

Note that Corollary~\ref{cor:Fk(n)extnotasp} generalizes Corollary~\ref{cor:F(n)extnotasp}.

\begin{corollary}\label{cor:Fk/l(3)extnotasp}
Suppose $k,l\geq 1$, $(k,l)=1$. Then the relative presentation $\mathcal{P}=\pres{G,x}{x^ltx^ktx^{-l}t^{-2}}$ (where $G=\pres{t}{t^3}$) is not aspherical.
\end{corollary}

\begin{proof}
By Corollary~\ref{cor:Fkl(n)extnotasp} we may assume $w(3,k)=1$, so $F^{k/l}(3)$ is abelian, by Corollary~\ref{cor:Fk/l(3)abelian}. But then $F^{k/l}(3)$ is finite, so the result follows from Corollary~\ref{cor:Fkl(n)extnotasp}.
\end{proof}

We now introduce the following quotients of groups $F^{k/l}(n)$. For each $n\geq 2$, $k,l\geq 1$ and each $\Omega\geq 0$ define
\[ F^{k/l}(n;\Omega)= \pres{x_0,\ldots , x_{n-1}}{x_i^lx_{i+1}^k=x_{i+2}^l, x_i^\Omega=1\ (0\leq i<n)}.\]
\begin{lemma}\label{lem:nonabelianfiniteordergens}
Let $m\geq 3$, $K,L\geq 1$, $(K,L)=1$, $\Omega\geq 0$. Suppose $F^{K/L}(m;\Omega)$ is infinite (resp.\,is non-cyclic, resp.\,is non-abelian, resp.\,is non-solvable). Then for all $n,k,l$ where $k\equiv \pm K$~mod~$\Omega$, $l\equiv \pm L$~mod~$\Omega$, $n\equiv 0$~mod~$m$, the group $F^{k/l}(n)$ is infinite (resp.\,is non-cyclic, resp.\,is non-abelian, resp.\,is non-solvable). Further, if $n/m$ is odd and $w(m,K)\neq 1$ in $F^{K/L}(m;\Omega)$ then $w(n,k)\neq 1$ in $F^{k/l}(n)$.
\end{lemma}

\begin{proof}
Let $\epsilon=\pm 1,\delta =\pm 1$, $n\equiv 0$~mod~$m$, $k\equiv \epsilon K$~mod~$\Omega$, $l\equiv \delta L$~mod~$\Omega$. Let $\phi:F^{k/l}(n) \rightarrow F^{k/l}(m)$ be the natural epimorphism given by $\phi(x_i)=x_{i~\mathrm{mod}~m}$. We have $F^{k/l}(m)\cong F^{\epsilon k/\delta l}(m)$ by Lemma~\ref{lem:F^k/lsigns}, so it maps onto $F^{K/L}(m;\Omega)$. If this latter group is infinite (or is non-cyclic or is non-abelian, or is non-solvable) then the same therefore holds for $F^{k/l}(n)$. It remains to show that if $n/m$ is odd and $w(m,K)\neq 1$ in $F^{K/L}(m;\Omega)$ then $w(n,k)\neq 1$ in $F^{k/l}(n)$.

Now if $n/m$ is odd then $\phi (w(n,k))=w(m,k)^{n/m}= (w(m,k)^2)^{(n/m-1)/2} w(m,k)=w(m,k)\in F^{k/l}(m)$. By adjoining the relators $x_i^\Omega$ ($0\leq i<m$) the group $F^{k/l}(m)$ maps onto $F^{\epsilon K/\delta L}(m;\Omega)\cong F^{K/L}(m;\Omega)$. Thus if $w(m,K)\neq 1$ in $F^{K/L}(m;\Omega)$ then $w(m,K)\neq 1$ in $F^{k/l}(m)$ and hence $w(n,k)\neq 1$ in $F^{k/l}(n)$.
\end{proof}

In Corollaries~\ref{cor:order2n3mod6}, \ref{cor:wneq1n3mod6}, \ref{cor:m=5or7infab} we give applications of Lemma~\ref{lem:nonabelianfiniteordergens} and in Example~\ref{ex:k/l(n;q)} we give further examples of groups $F^{k/l}(m;\Omega)$ to which Lemma~\ref{lem:nonabelianfiniteordergens} can usefully be applied.

\begin{corollary}\label{cor:order2n3mod6}
Suppose $n\equiv 3$~mod~$6$, $k,l\geq 1$, $(k,l)=1$. If $F^{k/l}(3)$ is non-abelian then $w(n,k)\neq 1$ in $F^{k/l}(n)$.
\end{corollary}

\begin{proof}
By Corollary~\ref{cor:Fk/l(3)abelian} if  $F^{k/l}(3)$ is non-abelian then $w(3,k)\neq 1$ in $F^{k/l}(3)=F^{k/l}(3;0)$ so the result follows from Lemma~\ref{lem:nonabelianfiniteordergens}.
\end{proof}

In cases where the order of the generators $x_i$ of $F^{K/L}(m)$ are known and finite we can set $\Omega$ equal to that order. However, it can be fruitful to set $\Omega$ to be a proper divisor of that order. Both instances are exhibited in the proof of the following corollary, where the order of generators $x_i$ of $F^{1/1}(3)$ is equal to 4 (and we set $\Omega=4$); whereas the order of the generators $x_i$ of the groups $F^{1/3}(3)$, $F^{1/4}(3)$, $F^{2}(3)$, $F^{3/2}(3)$ is 28, 49, 14, 63, respectively (and we set $\Omega=7$).

\begin{corollary}\label{cor:wneq1n3mod6}
\begin{itemize}
  \item[(a)] If $k$ and $l$ are odd and $n\equiv 3$~mod~$6$ then $w(n,k)\neq 1$ in $F^{k/l}(n)$.
  \item[(b)] If $(\pm k~\mathrm{mod}~7, \pm l~\mathrm{mod}~7)\in \{ (1,3), (2,1), (3,2)\}$ and $n\equiv 3$~mod~$6$ then $w(n,k)\neq 1$ in $F^{k/l}(n)$.
\end{itemize}
\end{corollary}

\begin{proof}
(a) This follows from Lemma~\ref{lem:nonabelianfiniteordergens} by observing that $F^{1/1}(3;4)\cong Q_8$ and so $w(3,1)\neq 1$ in this group by Corollary~\ref{cor:Fk/l(3)abelian}.
(b) This follows from Lemma~\ref{lem:nonabelianfiniteordergens} by observing that $F^{1/3}(3;7)\cong F^{2/1}(3;7) \cong F^{3/2}(3;7)$ is a non-abelian group (of order 56) and so $w(3,1)\neq 1$ in this group by Corollary~\ref{cor:Fk/l(3)abelian}.
\end{proof}

\begin{corollary}\label{cor:m=5or7infab}
If $(k,l)=1$, $k$ is even, $l\equiv 3$~mod~$6$, and $n\equiv 0$~mod~$m$, where $m\in \{5,7\}$ then $F^{k/l}(n)$ is infinite.
\end{corollary}

\begin{proof}
The hypotheses imply $k\equiv \pm 2$~mod~$6$ and $l\equiv 3$~mod~$6$. For $m\in \{5,7\}$ computations in GAP show that $F^{2/3}(m;6)$ has an index 5 subgroup with infinite abelianisation. Therefore $F^{2/3}(m;6)$ is infinite, and the result follows from Lemma~\ref{lem:nonabelianfiniteordergens}.
\end{proof}

\begin{example}\label{ex:k/l(n;q)}
\em
\begin{itemize}
  \item[(a)] $F^{1/3}(5;6)\cong PSL(2,11)$; $F^{3/1}(3;36)$ is a non-abelian, solvable group of order 3528; $F^{3/2}(3;63)$ is a non-abelian, solvable group of order 504; $F^{3/1}(3;6)\cong \Z_2\oplus \Z_6$.

  \item[(b)] Computations with the {\tt{NewmanInfinityCriterion}} command \cite{Newman} in GAP (applied to the derived subgroup or second derived subgroup) show that the following groups $F^{k/l}(n;\Omega)$ are infinite, and therefore $w(n,k)\neq 1$ in $F^{k/l}(n)$ by Corollaries~\ref{cor:Fkabelian} and~\ref{cor:Fk/l(3)abelian}: $F^{1/1}(9;76)$, $F^{3/1}(9;108)$, $F^{2/1}(17;206)$, $F^{2/1}(21;98)$, $F^{2/1}(23;94)$, $F^{3/5}(3;126)$, $F^{1/11}(3;182)$. The $\Omega$ values are selected as divisors of the order of $F^{k/l}(n)^\mathrm{ab}$ that are large enough for the quotient $F^{k/l}(n,\Omega)$ to be infinite yet small enough for the {\tt{NewmanInfinityCriterion}} command to complete.
\end{itemize}
\em
\end{example}

\section{3-manifold groups}\label{sec:3mfd}

Theorem~\ref{thm:Worder2inF(n)} was used in \cite{HowieWilliamsMFD} to obtain the following result.

\begin{theorem}[{\cite[Theorem~3]{HowieWilliamsMFD}}]\label{thm:F(n)not3mfd}
If $n\geq 3$ is odd then $F(n)$ is a 3-manifold group if and only if $n=3,5,7$, in which case $F(n)\cong Q_8, \Z_{11}, \Z_{29}$, respectively.
\end{theorem}

In this section we use the results of Section~\ref{sec:torsion} to prove the following corresponding result to Theorem~\ref{thm:F(n)not3mfd} for the groups $F^k(n)$ and $F^{k/l}(3)$. As reported earlier, the groups $F(5),F(7),F^{1/2}(3), F^{2/3}(3)$ are cyclic and $F(3)\cong Q_8$, and we expect this to be the only non-cyclic 3-manifold group among the groups $F^k(n)$ and $F^{k/l}(3)$. Part~(a) of Theorem~\ref{thm:Fkl(3)Fk(n)3mfd} is in contrast to the case when $n$ is even, where $F^{k}(n)$ is the fundamental group of a hyperbolic 3-manifold if either $k=1$ and $n\geq 8$ \cite[Theorem~C]{HKM} or $k\geq 2$ and $n\geq 6$ \cite[Theorem~3]{MaclachlanReid}.

\begin{theorem}\label{thm:Fkl(3)Fk(n)3mfd}
\begin{itemize}
  \item[(a)] Let $n\geq 3$ be odd, $k\geq 1$. If $F^k(n)$ is a non-cyclic 3-manifold group then $k$ is odd, $n\equiv 3$~mod~$6$ and $F^k(n)\cong Q_8\times \Z_{V^{k/1}_n/4}$.
  \item[(b)] Let $k,l\geq 1$, $(k,l)=1$. If $F^{k/l}(3)$ is a non-cyclic 3-manifold group then $k$ and $l$ are odd and $F^{k/l}(3)\cong Q_8\times \Z_{V^{k/l}_3/4}$.
\end{itemize}
\end{theorem}

We first extract an argument from the proof of Theorem~\ref{thm:F(n)not3mfd} and apply it to groups $F^{k/l}(n)$:

\begin{lemma}\label{lem:inf3mfd}
Let $n\geq 3$ be odd, $k,l\geq 1$, let $G=F^{k/l}(n)$ and let $w(n,k)=x_0^kx_1^k\ldots x_{n-1}^k\in G$. If $w(n,k)\neq 1$ then $G$ is not a 2-generator, infinite, 3-manifold group. In particular:
\begin{itemize}
  \item[(a)] if $n\geq 3$ is odd and $k\geq 1$ then $F^k(n)$ is not an infinite 3-manifold group;
  \item[(b)] if $k,l\geq 1$, where $(k,l)=1$, then $F^{k/l}(3)$ is not an infinite 3-manifold group.
\end{itemize}
\end{lemma}

\begin{proof}
Suppose $G$ is a 2-generator, infinite, 3-manifold group. By Theorem~\ref{thm:order2inFk/l} we have $w(n,k)\in G'$. Therefore the subgroup $<w(n,k)>\cong\Z_2$ is an orientation preserving subgroup of $G=\pi_1 (M)$  of finite order. Then by \cite[Theorem~8.2]{Epstein} (see also \cite[Theorem~9.8]{Hempel}) we have $M= R \# M_1$ where $R$ is closed and orientable, $\pi_1(R)$ is finite, and $<w(n,k)>$ is conjugate to a subgroup of $\pi_1(R)$. Since $G$ is infinite we have $\pi_1(M_1)\neq 1$ and since it can be generated by two elements $\pi_1(R)$ and $\pi_1(M_1)$ are each cyclic. But the derived subgroup of a free product of cyclic groups is free, contradicting the fact that $w(n,k)\in G'$ is an element of order two.

Part (a) (resp.\,Part~(b)) follows since $F^{k}(n)$ (resp.\,$F^{k/l}(3)$) is 2-generated by Lemma~\ref{lem:Fk2gen} (resp. Lemma~\ref{lem:(k,l)=1isn+1/2gen}) and if it is infinite then $w(n,k)\neq 1$ by Theorem~\ref{thm:Fk/labVn} and Corollary~\ref{cor:Fkabelian} (resp.\,Corollary~\ref{cor:Fk/l(3)abelian}).
\end{proof}

To consider when $F^{k}(n)$ and $F^{k/l}(3)$ can be finite 3-manifold groups we need the following classification of finite 3-manifold groups (see \cite[Section~2]{Hopf} or \cite[Section~1.5]{Wilton}) and their derived subgroups.

\begin{theorem}\label{thm:finite3mfdgroups}
Suppose $G$ is a finite 3-manifold group. Then either $G$ is cyclic or $G\cong H\times \Z_p$ where $p\geq 1$ is coprime to $|H|$ and $H$ is as in one of the following cases:
\begin{itemize}
  \item[(i)] $H=P_{48}=\pres{x,y}{x^2=(xy)^3=y^4, x^4=1}$, with $H/H'\cong \Z_2$, $H'\cong SL(2,3)$ and $H'/H''\cong \Z_3$;
  \item[(ii)] $H=P_{120}=\pres{x,y}{x^2=(xy)^3=y^5, x^4=1}$, a perfect group;
  \item[(iii)] $H=Q_{4m}=\pres{x,y}{x^2=(xy)^2=y^m}$, $m\geq 2$, with $H/H'\cong \Z_2\oplus \Z_2$ ($m$ even), $H/H'\cong \Z_4$ ($m$ odd), and $H'\cong \Z_m$;
  \item[(iv)] $H=D_{2^m(2n+1)}=\pres{x,y}{x^{2^m}=1, y^{2n+1}=1, xyx^{-1}=y^{-1}}$, $m,n\geq 1$, with $H/H'\cong \Z_{2^m}$ and $H'\cong \Z_{2n+1}$;
  \item[(v)] $H=P'_{8\cdot 3^m}=\pres{x,y,z}{x^{2}=(xy)^2=y^2, zxz^{-1}=y, zyz^{-1}=xy, z^{3^m}=1}$, $m\geq 1$, with $H^\mathrm{ab}\cong \Z_{3^m}$ and $H'\cong Q_8$.
  \end{itemize}
\end{theorem}

We now prove Theorem~\ref{thm:Fkl(3)Fk(n)3mfd}.

\begin{proof}[Proof of Theorem~\ref{thm:Fkl(3)Fk(n)3mfd}]
Let $G=F^{k/l}(n)$ where $k,l\geq 1$, $(k,l)=1$ and either $l=1$ or $n=3$. By Lemma~\ref{lem:inf3mfd} we may assume that $G$ is a finite, non-cyclic, 3-manifold group.

Observe that in each case $w=w(n,k)\neq 1$ by Corollaries~\ref{cor:Fkabelian} and~\ref{cor:Fk/l(3)abelian}, and that $G$ is generated by $x_0$ and $x_1$ by Lemmas~\ref{lem:Fk2gen} and~\ref{lem:(k,l)=1isn+1/2gen}. Suppose that $G \cong H\times \Z_p$ where $H$ is one of the groups in (i)--(v) of Theorem~\ref{thm:finite3mfdgroups} and $(p,|H|)=1$. Then the derived subgroup $D$ of $G$ is isomorphic to the derived subgroup of $H$. By Theorem~\ref{thm:order2inFk/l} the element $w(n,k)=[x_0^l,x_1^l]$ has order 2 in $G$. Corollaries~\ref{cor:Fkabelian} and~\ref{cor:Fk/l(3)abelian} imply that $D$ is the normal closure of $w$ in $G$ and so $D^\mathrm{ab}\cong \Z_2^d$ for some $d\geq 0$, which gives a contradiction if $H$ is the group in part (i) or (iv).

If $H=P_{120}$ or  $H=P'_{8\cdot 3^m} \cong Q_8 \rtimes \Z_{3^m}$ (as in parts (ii) or (v)) then $H$ has a unique element $h$ of order 2, and the normal closure $\npres{h}^H$ is not isomorphic to the derived subgroup of $H$. Therefore $H\times \Z_p$ has a unique element of order 2, namely $(h,0)$, and the normal closure $\npres{(h,0)}^{H\times \Z_p}$ is not the derived subgroup of $H\times \Z_p\cong G$, a contradiction (since the normal closure $\npres{w}^G=D$).

If $H=Q_{4m}$ for some $m\geq 2$ (as in part~(iii)), then $D\cong \Z_m$, so $m=2$ and hence $G\cong Q_8\times \Z_p$, where $p=V_n^{k/l}/4$.
Therefore there is an epimorphism $G\twoheadrightarrow\Z_2 \oplus \Z_2$, and since the images of generators $x_0,x_1$ have order 2 there is also an epimorphism $F^{k/l}(n;2) \twoheadrightarrow \Z_2\oplus \Z_2$. If $k$ is even then $(k,l)=1$ implies $l$ is odd so $F^{k/l}(n;2)\cong \Z_2$, a contradiction. Therefore $k$ is odd. Since $|G^\mathrm{ab}|=V_n^{k/l}$ is even the remaining conditions on $k,l,n$ follow from Corollary~\ref{cor:Vneven}.
\end{proof}

\section*{Acknowledgements}

The authors thank Jim Howie for helpful comments on a draft and the referee for the careful reading of the article.

  \textsc{Department of Mathematical Sciences, University of Essex, Wivenhoe Park, Colchester, Essex CO4 3SQ, UK.}\par\nopagebreak
  \textit{E-mail address}, \texttt{Ihechukwu.Chinyere@essex.ac.uk}

  \textsc{Department of Mathematical Sciences, University of Essex, Wivenhoe Park, Colchester, Essex CO4 3SQ, UK.}\par\nopagebreak
  \textit{E-mail address}, \texttt{Gerald.Williams@essex.ac.uk}

\end{document}